\documentclass[11pt, reqno]{amsart}
\usepackage{amsmath, amsthm, a4, latexsym, amssymb}

\def\@rmrk#1#2{\refstepcounter
    {#1}\@ifnextchar[{\@yrmrk{#1}{#2}}{\@xrmrk{#1}{#2}}}

\makeatletter\@addtoreset{equation}{section}\makeatother

\newfont{\bfit}{cmbxti10 scaled 2000}
\newfont{\biggi}{cmr12 scaled 2000}
\newtheorem{step}{STEP}

\newcommand{\bes}{\begin{step}}
\newcommand{\es}{\end{step}}
 
 \newcommand{\eps}{\varepsilon}

 \newcommand{\R}{\mathbb{R}}
 
 \newcommand{\N}{\mathbb{N}}

 \newcommand{\prob}{\mathbb{P}}

 \newcommand{\skrib}{{\mathcal B}}
 \newcommand{\skric}{{\mathcal C}}

 \newcommand{\skrie}{{\mathcal E}}
 
 \newcommand{\skrig}{{\mathcal G}}

 \newcommand{\skril}{{\mathcal L}}
 \newcommand{\skrim}{{\mathcal M}}

 \newcommand{\skrin}{{\mathcal N}}
 \newcommand{\skrip}{{\mathcal P}}

 \newcommand{\sfrac}[2]{\mbox{$\frac{#1}{#2}$}}

\def\1{{\mathchoice {1\mskip-4mu\mathrm l}      
{1\mskip-4mu\mathrm l}
{1\mskip-4.5mu\mathrm l} {1\mskip-5mu\mathrm l}}}

\newcommand{\eq}{\begin{equation}}
\newcommand{\en}{\end{equation}}

\renewcommand{\subsection}{\secdef \subsct\sbsect}
\newcommand{\subsct}[2][default]{\refstepcounter{subsection}
\vspace{0.15cm}
{\flushleft\bf \arabic{section}.\arabic{subsection}~\bf #1  }
\nopagebreak\nopagebreak}
\newcommand{\sbsect}[1]{\vspace{0.1cm}\noindent
{\bf #1}\vspace{0.1cm}}

\newtheorem{theorem}{Theorem}[section]
\newtheorem{lemma}[theorem]{Lemma}
\newtheorem{cor}[theorem]{Corollary}

\newtheoremstyle{thm}{1.5ex}{1.5ex}{\itshape\rmfamily}{}
{\bfseries\rmfamily}{}{2ex}{}

\newtheoremstyle{rem}{1.3ex}{1.3ex}{\rmfamily}{}
{\itshape\rmfamily}{}{1.5ex}{}
\theoremstyle{rem}
\newtheorem{remark}{{\slshape\sffamily Remark}}[]

\refstepcounter{subsection}

\def\thebibliography#1{\section*{reference}
  \list%
  {\arabic{enumi}.}
    {\settowidth\labelwidth{[#1]}\leftmargin\labelwidth
    \advance\leftmargin\labelsep
    \parsep0pt\itemsep0pt
    \usecounter{enumi}}
    \def\newblock{\hskip .11em plus .33em minus .07em}
    \sloppy                   
    \sfcode`\.=1000\relax}

\begin{document}
\title[Large  deviation principle for empirical  measures of  Multitype random networks]
{\Large Large  deviation principle for empirical  measures of  Multitype random networks}

\author[K. Doku-Amponsah]{}

\maketitle
\thispagestyle{empty}
\vspace{-0.5cm}

\centerline{\sc{By Kwabena Doku-Amponsah}}
\renewcommand{\thefootnote}{}
\renewcommand{\thefootnote}{1}
\renewcommand{\thefootnote}{}
\footnote{\textit{Address:} Statistics Department, University of
Ghana, Box LG 115, Legon,Ghana.\,
\textit{E-mail:\,kdoku@ug.edu.gh}.}
\renewcommand{\thefootnote}{1}
\centerline{\textit{University of Ghana}}

\begin{quote}{\small }{\bf Abstract.} In  this  article  we  study  the stochastic  block  model also known  as  the  multi-type  random  networks (MRNs). For the stochastic block model or the MRNs  we define  the  empirical group  measure,  empirical cooperative  measure  and  the  empirical  locality  measure. We derive large  deviation  principles  for  the  empirical  measures in  the  weak  topology. These results will form the  basis of  understanding  asymptotics  of  the  evolutionary  and  co-evolutionary   processes on  the  stochastic block  model.

\end{quote}\vspace{0.5cm}

\section{Introduction}

Networks  are  complex  objects  and  may  be  sometimes  difficult  to  study.
The  role  of network structure in  information  Science,Social  Science  and behavioural  science is  well  understood. In  behavioural Science (such  as  Sociology, Economics) and  applied  statistics the study  of social  bonds among actors is a classical field,  known as  social  network analysis.  See Wasserman and Faust(1994). Modelling  the  evolution  of  networks  is  central  to our  understanding  of  large  communication systems,  and  more  general, modern economics and social  systems. See Hellman and  Staudigl (2014).   The  research  on  social  network  and  economic  networks  is  interdisciplinary and  the  number  of  proposed models is  huge.\\

 Often  the  sites  in  a  network  can  be  classified  to  belong  to  certain groups. A  most  common  feature  of  social  networks    is  the  phenomenom  of  homophily.  Thus,  sites  of  similar  attributes  are  more  likely  to  be  linked  with  each  other. Fienberg,  Meyer  and  Wasserman(1985)  introduced  stochastic blocking  models,  where   sites  are categorized  to  come  from  a   certain  subgroup. Every  subgroup may  have  it  own  law  of  network  formation. More  recently, this  type  of  networks  has  also  been  used  in  the  field  of  Economics,  see  example (Golub,  2012),  where  it  has  been  called  a  multitype  random networks.\\

In  this  paper  we  use  large  deviation  technique  to  obtain asymptotic results for  multi-type  random  networks. See, Doku-Amponsah and Moerters (2010), Doku-Anponsah (2015) for similar  results  for coloured random graphs.  To  be  specific, we find  in this  article   joint  large  deviation  principles  for  empirical group measure  and  the  empirical  co-operative  measure, and  the  empirical co-operative  measure  and  the  empirical locality  measure  of  the  multitype  random networks  in  the  weak  topology. We  obtain  these results  using  (Doku-Amponsah and Moerters 2010,  Theorem 2.1 and  Theorem 2.3).\\

 Our  main motivation  for  studying  this  model  is two-fold.  In  one  hand  we  try to  improve  the  presentations  of  the  LDPs  for  the  inhomogeneous  random graphs  by  other  researcher   and  on  the  other  hand  we  working  towards  understanding  the  asymptotics  of  the  evolution and  co-evolution  processes  on  random  networks.  Thus, LDPs  develop   for empirical  measures of the  multitype  random  networks  will form  the  basis for understanding  evolution  as  the  building block of the  construction  of  co-evolution model  of  networks  and  play.\\


 \subsection{Empirical  Measures  of  the Multitype Random  Networks}

 By  $\Omega$  we  denote a  finite  alphabet  and  denote   by $\skrin(\Omega)$  the  space  of  counting  measure  on   $\Omega$
 equipped  with  the  discrete topology. By  $\skrim(\Omega)$  we  denote  the  space of  probability  measures  on  $\Omega$  equipped  with  the weak  topology  and  $\skrim_*(\Omega)$  the  space
 of finite  measures  on  $\Omega$  equipped  with  the  weak  topology.\\

 Let  $\kappa_n,\ell_n:\Omega\times\Omega\to[0,\,\infty)$    be    continuous  functions and  $\eta:\Omega\to (0,\,1)$ a probability  measure.  The  Multitype Random  Network  is  defined  as  follows:
 \begin{itemize}
 \item[(i)] Partition the  set  of  sites, $[n]=\{1,2,3,...,n\},$ into  finitely  many  types  $\Omega=\{a_1,a_2,a_3,...,a_m\}$ independently  according  to  the  type  law  $\eta.$
 \item[(ii)] For  a  couple  of  sites $(i,j)$   take  the  intensities  of  link  formation  and link  destruction  to  be  $w_{ij}^{(n)}(a,b)=(1-A_{ij})\kappa_n(a,b),\,\, \nu_{ij}^{(n)}(a,b)=A_{ij}\ell_{n}(a,b)$, whenever site $i$  is  a  member  of  group $a$  , and   site  $j$  is  a  member  of  group  $b.$ 
\end{itemize}

We  shall  consider  $Z:=\{ (Z(i),\in \Omega),E\}$  under  the  joint  law  of  the  type  and  graph,   and  interpret  $Z(i)$  as  the  type  of  site  $i$  and $Z$  as  multi-type  random  graph.  We  shall  refer  to  the  multitype  random graph  $Z$   \emph{ as  symmetric } if  both  functions  $\kappa_n$  and  $\ell_n$  are  symmetric  otherwise  we  call it \emph{ asymmetric}  multitype  random  graph.

This  model is  a  special  case  of  the  inhomogeneous  random  graph,  see example  Hellman and  Staudigl \cite{HS14},  where  edge  specific   probabilities between  members  of  group  $a$  and  group  $b$  are given  by  $$p_n(a,b)=\sfrac{\kappa_{n}(a,b)}{\kappa_{n}(a,b)+\ell_{n}(a,b)}.$$

The  multi-type  random  graph,  unlike  the  inhomogeneous  random  graph, is  completely  specified  by  the grouping  of  the   sites, and  the  group-specific  link-success  probabilities.

With each  multitype random  graph $Z$  we  associate  a probability distribution, the
\emph{empirical type measure}~$L^1\in\skrim(\Omega)$,~by
$$L^{1}(a):=\frac{1}{n}\sum_{j=1}^{n}\delta_{Z(j)}(a),\quad\mbox{ for $a\in\Omega$, }$$
and a finite measure, the \emph{empirical cooperative measure}
$L^{2}\in\tilde\skrim_*(\Omega\times\Omega),$ by
$$\begin{aligned}
L^{2}(a,b)=\left\{
\begin{array}{ll}\frac{1}{n}\sum_{(i,j)\in E}[\delta_{(Z(i),Z(j))}+
\delta_{((Z(j),Z(i))}](a,b)& \mbox {if $Z$  is  symmtric }\\
\frac{2}{n}\sum_{(i,j)\in E}\delta_{(Z(i),Z(j))}(a,b) & \mbox{if $Z$ is  asymmetric.}
\end{array}\right.
\end{aligned}$$

Notice,$\|L^2\|=2|E|/n$  in  both the  symmetric  and  asymmetric cases  of  the  multitype  random  graph.  We define a probability distribution, the
\emph{empirical neighbourhood measure}
$M^1\in\skrim(\Omega\times\skrin(\Omega))$, by
$$M^1(a,l):=\frac{1}{n}\sum_{j=1}^{n}\delta_{(Z(j),L(j))}(a,l),\quad
\mbox{ for $(a,l)\in \Omega\times\skrin(\Omega)$, }$$ where
$L(i)=(l^{i}(b),\,b\in\Omega)$ and $l^{i}(b)$ is the number of
sites  of type $b$  linked to site $i$.	

  \section{Main  Results}\label{Main}
We  assume  through out the  paper  that  the  intensities of  link  formation  and  link  destruction satisfies 
 $n\kappa_{n}(a,b)/ \ell_{n}(a,b)\to \kappa(a,b)/\ell(a,b),$  while $\kappa,\ell:\Omega\times\Omega\to[0,\,\infty)$. 
We  write   $\sfrac{\kappa}{\ell}\rho\otimes\rho(a,b)=\sfrac{\kappa(a,b)}{\ell(a,b)}\rho(a)\rho(b)$ for
$a,b\in\Omega$  and  $${\mathfrak H}_{\kappa/\ell}(\pi\, \| \, \rho ):=
H\big(\omega\,\|\,\sfrac{\kappa}{\ell}\rho\otimes\rho\big)+\|
\sfrac{\kappa}{\ell}\rho\otimes\rho \| -\|\pi\|.\, $$
It is not hard to see that $\mathfrak
H_1(\pi\,\|\,\omega_1)\ge 0$ and equality holds if and only if
$\pi= \sfrac{\kappa}{\ell}\omega_1\otimes\omega_1$.



\begin{theorem}\label{LDP1}
Suppose  $G$  is  a  multitype  random  network  with grouping law $\eta:\Omega\to[0,1]$. Assume,that  the  intensities of  link  formation  and  link  destruction  are $$w_{ij}(a,b)=(1-A_{ij})\kappa_n(a,b),\,\,\, \nu_{i,j}(a,b)=A_{ij}\ell_n(a,b),$$   where $\kappa:\Omega\times\Omega\to [0,\,\infty)$ and $\ell:\Omega\times\Omega\to(0,\infty],$  whenever site $i$  is  a member  of  group  $a$  and  site  $j$  is  a  member  of  group $b,$  for  $a,b\in \Omega.$  Then,  the  pair  $(L^2, M)$  an  LDP  obeys in  the  space  $\skrim(\Omega\times\Omega)\times\skrim(\Omega\times\skrin(\Omega))$  with  speed  $n$  and  with good rate function
$$\begin{aligned}
J^1(\pi, \omega)=\left\{
\begin{array}{ll}H(\omega\,\|\,q^{1})+H(\omega_1\,\|\,\eta)+\sfrac{1}{2}{\mathfrak H}_2(\pi\|\omega_1)& \mbox {if $(\pi,\omega)$ consistent  and  $\omega_1=\pi_2,$   }\\
\infty & \mbox{otherwise.}
\end{array}\right.
\end{aligned}$$
where  $$q^1(a,\rho):=\omega_1(a)\prod_{b\in\Omega}\frac{e^{-\big[\pi(a,b)/\mu_1(a)\big]}{\big[\pi(a,b)/\omega_1(a)\big]}^{\rho(b)}}{\rho(b)!}$$



\end{theorem}

\begin{cor}\label{LDP2}
	Suppose  $deg_Z$ is  the  degree  measure of  $Z$,  an  Erdos-Renyi  graph model. Assume,that  the  intensities of  link  formation  and  link  destruction  are $$w_{ij}=(1-A_{ij})\kappa,\,\,\, \nu_{i,j}=A_{ij}\ell.$$ 
	\begin{itemize}
		
		\item	Then , $deg_Z$
		satisfies an LDP in the space $\skrim(\N \cup \{0\})$ with   speed $n$   and  good rate
		function
		\begin{equation}\label{randomg.ratedeg}
		\begin{aligned}
		\lambda(d)= \left\{ \begin{array}{ll}\Big[ H (d\,\|\,q_{\langle d\rangle})+\sfrac 12 \, \langle
		d \rangle\, \log \big( \sfrac{\ell\langle d\rangle}{\kappa}
		\big)
		- \sfrac 12 \, \langle d \rangle  +  \sfrac {\kappa}{2\ell}\Big] ,
		& \mbox { if $\langle d\rangle< \infty,$ }\\[2mm]
		\infty & \mbox { if $\langle d\rangle= \infty.$ }
		\end{array}\right.
		\end{aligned}
		\end{equation}
		where $q_{c}(k)=\sfrac{e^{-c}c^k}{k!},$  $k=0,1,2,...$ and  $\langle d \rangle =\sum_{k=0}^{\infty}kd(k).$
		\item  Then,  the  proportion of isolated vertices, $deg_Z(0)$ satisfies an LDP in $[0,1]$
		with good rate function $$h(z)=z \log z + \sfrac{\kappa}{\ell} z(1-z/2) -
		(1-z) \big[ \log\big (\sfrac {\kappa}{\ell a}\big)-
		\sfrac{(a-  \sfrac{\kappa}{\ell} (1-z))^2}{2\sfrac{\kappa}{\ell}(1-z)} \big] \, ,$$
		where $t=t(z)$ is the unique positive solution of $1-e^{-t}=\frac {\kappa}{\ell t}\, (1-z)$.
		
	\end{itemize}
\end{cor}

\begin{remark}
	Note,  typically,  as  $n\to\infty$,  the  proportion  of  isolated  sites  in  an  Erdos-Renyi  graph      converges  to  $e^{-n\kappa/\ell}$  in probability.i.e. $$\lim_{n\to\infty}\prob\Big\{|deg_X(0)-e^{-n\kappa/\ell}|\ge \eps\Big\}=0.$$
\end{remark}
\begin{remark}
	Note,  typically,  as  $n\to\infty$,  the  proportion  of  isolated  sites  in  an  Erdos-Renyi  graph      converges  to  $e^{-n\kappa/\ell}$  in probability.i.e. $$\lim_{n\to\infty}\prob\Big\{|deg_X(0)-e^{-n\kappa/\ell}|\ge \eps\Big\}=0.$$
\end{remark}

\begin{theorem}\label{LDP4}
	Suppose  $Z$  is  a  multitype  random  network  with grouping law $\eta:\Omega\to[0,1]$. Assume,that  the  intensities of  link  formation  and  link  destruction  are $$w_{ij}(a,b)=(1-A_{ij})\kappa(a,b),\,\,\, \ell_{i,j}(a,b)=A_{ij}\ell(a,b),$$   where $\kappa:\Omega\times\Omega\to [0,\,\infty)$ and $\ell:\Omega\times\Omega\to(0,\infty],$  whenever site $i$  is  a member  of  group  $a$  and  site  $j$  is  a  member  of  group $b$  for  $a,b\in \Omega.$ 
	Then,  $L^2$  obeys  the  following   large  deviation  principle:
	\begin{itemize}
		\item  for  any  $F$ closed  subset of $\skrim(\Omega\times\Omega)$,  we have  
		$$\limsup_{n\to \infty}\frac{1}{n}\log P_{\omega}\Big\{L^2\in F\Big\}\le -\inf_{\pi\in F}I_{1}(\omega, \pi)$$
		\item for  any  $\Gamma$ open  subset  of $\skrim(\Omega\times\Omega)$, we  have    
		$$\liminf_{n\to \infty}\frac{1}{n}\log P_{\omega}\Big\{L^2\in F\Big\}\ge -\inf_{\pi\in F}I_{1}(\omega, \pi)$$
	\end{itemize}
\end{theorem} 

	
		

\begin{remark}
	Note,  typically,  as  $n\to\infty$,  the  proportion  of  isolated  sites  in  an  Erdos-Renyi  graph      converges  to  $e^{-n\kappa/\ell}$  in probability.i.e. $$\lim_{n\to\infty}\prob\Big\{|deg_X(0)-e^{-n\kappa/\ell}|\ge \eps\Big\}=0.$$
\end{remark}

\begin{remark}
We  observe  that  as  $n\to\infty$,  given  $\langle d\rangle=\kappa/\ell$  the  degree  distribution converges  $q_{\kappa/\ell}$ in  probability. i.e. $$\lim_{n\to\infty}\sup_{k\in\N\cap \{0\}}\prob\Big\{\Big|deg_Z(k)-\sfrac{e^{-\kappa/\ell}(\kappa/\ell)^{k}}{k!}\Big |\ge  \eps\Big\}=0.$$
\end{remark}

\begin{theorem}\label{LDP3}
Suppose  $Z$  is  a  multitype  random  network  with grouping law $\eta:\Omega\to[0,1]$. Assume,that  the  intensities of  link  formation  and  link  destruction  are $$w_{ij}(a,b)=(1-A_{ij})\kappa(a,b),\,\,\, \ell_{i,j}(a,b)=A_{ij}\ell(a,b),$$   where $\kappa:\Omega\times\Omega\to [0,\,\infty)$ and $\ell:\Omega\times\Omega\to(0,\infty],$  whenever site $i$  is  a member  of  group  $a$  and  site  $j$  is  a  member  of  group $b$  for  $a,b\in \Omega.$  Then,  the  pair  $(L^1,\,L^2)$  obeys an  LDP   in  the  space  $\skrim(\Omega)\times\skrim(\Omega\times\Omega)$  with  speed  $n$  and with good rate function   $$I(\rho,\,\pi)=H(\rho\,\|\,\eta)+\sfrac{1}{2}{\mathfrak H}_{\kappa/\ell}(\pi\, \| \, \rho )$$
\end{theorem}

\begin{remark}
We  observe  that  as  $n\to\infty$, given  $L^1=\rho$,  the  number  edges  per  site  $L^2$  converges  $\sfrac{\kappa}{\ell}\rho\otimes \rho$ in  probability. i.e. $$\lim_{n\to\infty}\prob\Big\{\Big|L^2-\sfrac{\kappa}{\ell}\rho\otimes \rho\Big |\ge  \eps\Big| L^1=\rho\Big\}=0.$$
\end{remark}
We write  $I_1(\rho,\,\pi)=\sfrac{1}{2}{\mathfrak H}_{\kappa/\ell}(\pi\, \| \, \rho )$ and  state  the following  Large  deviation  principle  for  $L^2.$

\begin{theorem}\label{LDP4}
	Suppose  $Z$  is  a  multitype  random  network  with grouping law $\eta:\Omega\to[0,1]$. Assume,that  the  intensities of  link  formation  and  link  destruction  are $$w_{ij}(a,b)=(1-A_{ij})\kappa(a,b),\,\,\, \ell_{i,j}(a,b)=A_{ij}\ell(a,b),$$   where $\kappa:\Omega\times\Omega\to [0,\,\infty)$ and $\ell:\Omega\times\Omega\to(0,\infty],$  whenever site $i$  is  a member  of  group  $a$  and  site  $j$  is  a  member  of  group $b$  for  $a,b\in \Omega.$ 
	Then,  $L^2$  obeys  the  following   large  deviation  principle:
	\begin{itemize}
		\item  for  any  $F$ closed  subset of $\skrim(\Omega\times\Omega)$,  we have  
		$$\limsup_{n\to \infty}\frac{1}{n}\log P_{\omega}\Big\{L^2\in F\Big\}\le -\inf_{\pi\in F}I_{1}(\omega, \pi)$$
		\item for  any  $\Gamma$ open  subset  of $\skrim(\Omega\times\Omega)$, we  have    
		$$\liminf_{n\to \infty}\frac{1}{n}\log P_{\omega}\Big\{L^2\in F\Big\}\ge -\inf_{\pi\in F}I_{1}(\omega, \pi)$$
	\end{itemize}
	\end{theorem}

\section{Proof  of  Theorem~\ref{LDP4}}
The  technique  use  in  this  section  is  routed  in  spectral  potential  theory  and  it is  a  summary  of  many  large  deviation  theories. To  begin,  for  any  multitype  random  network  $z$  conditioned  to  have  empirical  type  measure  $L^1=\omega$ we  define  the  spectral  potential 
$$\rho_{\kappa/\ell}(g,\omega)=-\langle (1-e^{g}),\sfrac{\kappa}{\ell}\omega\otimes\omega\rangle$$

We  observe  from  \cite{DA17d}  that  $\rho_{\kappa/\ell}$  possesses  all  the  remarkable  properties mention  by \cite{BIV15}  and  the  following  Lemma  holds  for  $L^2.$
\begin{lemma}\label{LDP5}
Suppose  $Z$  is  a  multitype  random  network  with grouping law $\eta:\Omega\to[0,1]$. Assume,that  the  intensities of  link  formation  and  link  destruction  are $$w_{ij}(a,b)=(1-A_{ij})\kappa(a,b),\,\,\, \ell_{i,j}(a,b)=A_{ij}\ell(a,b),$$   where $\kappa:\Omega\times\Omega\to [0,\,\infty)$ and $\ell:\Omega\times\Omega\to(0,\infty],$  whenever site $i$  is  a member  of  group  $a$  and  site  $j$  is  a  member  of  group $b$  for  $a,b\in \Omega.$ 
 Then,  $L^2$  obeys  the  following  local large  deviation  principle:
\begin{itemize}
	\item  for  any  $\pi\in\skrim(\Omega\times\Omega)$  and  a  number  $\eps>0$,  there  is a   weak  neighbourhood  $B_{\pi}$  such  that
	$$P_{\omega}\Big\{L^2\in B_{\pi}\Big\}\le e^{-nI_{1}(\omega, \pi)-n\eps+o(n)}$$
	\item for  any  $\pi\in\skrim(\Omega\times \Omega)$, any  number  $\eps>0$ and  a  fine  neighbourhood  $B_{\pi}$,  we  have  the  asymptotic  estimate 
	 $$P_{\omega}\Big\{L^2\in B_{\pi}\Big\}\ge e^{-nI_{1}(\omega, \pi)+n\eps-o(n)}$$
	\end{itemize}
\end{lemma}

Lemma~\ref{Lem1}  below,  which provides  a summaries  of  the   properties  of  the  Kullback  action  was  first  presented  in  the  paper \cite{DA17d}. To  state it, we denote by  $\skric$   the  space  of  continuous  functions    $g:\Gamma\times\Gamma\to \R.$

\begin{lemma}[\cite{DA17d}]\label{Lem1}\label{LLDP1}\label{LDP6}.The  following  holds  for  the  Kullback action or divergence  function  ${\mathfrak  H}_{\kappa/\ell}(\pi\,\|\,\omega).$
	\begin{itemize}
		\item [(i)]$\displaystyle {\mathfrak  H}_{\kappa/\ell}(\pi\,\|\,\omega)=\sfrac{1}{2}\sup_{g\in\skric}\Big\{\langle g,\,\pi\rangle-\rho_{\kappa/\ell}(g,\,\omega) \Big\}.$
		\item[(ii)] The  function ${\mathfrak  H}_{\kappa/\ell}(\pi\,\|\,\omega)$  is lower  semi-continuous on  the  space  $\skrip_*(\Gamma\times\Gamma).$
		\item[(iii)] For  any  real  $c,$  the  set  $\Big\{\nu\in\skrip_*(\Gamma\times\Omega):\, {\mathfrak  H}_{\kappa/\ell}(\pi\,\|\,\omega)\le  c\Big\}$  is  weakly  compact.
	\end{itemize}
\end{lemma}
Please  we refer to \cite{BIV15}  for  similar  result  and  proof  for  the  empirical  measures on  measurable spaces.


Note that  Lemma~~\ref{Lem1} (i) above  implies  the  so-called  variational  principle.See,  example ~\cite{KV86}.

\subsection{Proof  of  Lemma~\ref{LDP5}.}

We  write  $p_n(a,b):=\sfrac{\kappa_n(a,b)}{\kappa_n(a,b)+\ell_n(a,b)}$  and  $\tilde{p}_n(a,b):=\sfrac{\tilde{\kappa}_n(a,b)}{\tilde{\kappa}_n(a,b)+\tilde{\ell}_n(a,b)}.$  Recall  the  definition  of  $\tilde{h}_n $  from  \cite{DA12},and 
note  from  Lemma~\ref{LLDP1} that,  for  any  $\eps>0$  there  exists  a  function  $g\in\skrim(\Omega\times\Omega)$  such  that  ${\mathfrak  H}_{\kappa/\ell}(\pi\,\|\,\omega)-\sfrac{\eps}{2} < \langle g,\,\pi\rangle-\rho_{\kappa/\ell}(g,\,\omega) .$  We  define  the  probability  distribution $\tilde{P}_n$   by

\begin{align}
\tilde{P}_{\omega}(Z)
&= \prod_{(i,j)\in E}\tilde{p}_{n}(Z(i),Z(j))\prod_{(i,j)\not\in
	E}{1-\tilde{p}_{n}(Z(i),Z(j))}\nonumber\\
& = \prod_{(i,j)\in E}\sfrac{\tilde{p}_{n}(Z(i),Z(j))}{n-n\tilde{p}_{n}(Z(i),Z(j))}\prod_{(i,j)\in\skrie}{(n-n\tilde{p}_{n}(Z(i),Z(j)))}\nonumber\\
& = \prod_{(i,j)\in
	E}e^{g(Z(i),Z(j))}\sfrac{\tilde{p}_{n}(Z(i),Z(j))}{n-n\tilde{p}_{n}(Z(i),Z(j))}\prod_{(i,j)\in\skrie}{e^{\frac 1n\,
		\tilde{h}_n(Z(i),Z(j))}}{(n-n\tilde{p}_{n}(Z(i),Z(j)))}\label{changem1}
\end{align}

Using  Equation \ref{changem1}  above,  we  have

\begin{equation}
\frac{dP_{\omega}(Z)}{d\tilde{P}_{\omega}(Z)}
= \prod_{(i,j)\in
	E}e^{-g(Z(i),Z(j))}\prod_{(i,j)\in\skrie}{e^{-\frac 1n\,
		\tilde{h}_n(Z(i),Z(j))}}
= e^{-n\langle
	\sfrac{1}{2}L^2,\, \tilde{g}\rangle-n\langle\sfrac{1}{2}L^1\otimes
	L^1,\, \tilde{h}_n\rangle+\langle \sfrac{1}{2}L_{\Delta}^{1},\,
	\tilde{h}_n\rangle}/2,\label{changem2}
\end{equation}
while  $$L_{\Delta}^{1}=\delta_{(Z(i), Z(i))}.$$

Now  we  define a  neighbourhood  of  the  functional  $\pi$  as  follows:
$$B_{\pi}=\Big\{\varpi\in\skrim(\Omega\times\Omega):  \langle g, \, \varpi\rangle>\langle g, \, \pi\rangle -\sfrac{\eps}{2}\Big\}.$$
Therefore,  under  the  condition  $L^1\in B_{\pi}$    we  have  that  $$ \frac{dP_{\omega}(z)}{d\tilde{P}_{\omega}(z)}<e^{\sfrac{1}{2}(\rho_{\kappa/\ell}(g,\,\omega)-\langle g, \, \pi\rangle) +\sfrac{\eps}{2}}<e^{-n{\mathfrak  H}_{\kappa/\ell}(\pi\,\|\,\omega)+n\eps}.$$

Hence,  we have

$$\begin{aligned}
P_\omega\Big\{z\in\skrig([n],\,\Omega)| L_z^2\in B_{\pi}\Big\}\le \int_{\skrig([n],\Omega)}\1_{\{L_z^2\in B_{\nu}\}} d\tilde{P}_{\omega}(z) \le \int_{\skrig([n],\,\Omega)}\1_{\{L_y^2\in B_{\pi}\}} &e^{-n{\mathfrak  H}_{\kappa/\ell}(\pi\,\|\,\omega)-n\eps}d\tilde{P}_{\omega}(z) \\
&\le e^{-n{\mathfrak  H}_{\kappa/\ell}(\pi\,\|\,\omega)/2-n\eps} .
\end{aligned}$$
Note  that  $ {\mathfrak  H}_{\kappa/\ell}(\pi\,\|\,\omega)=\infty$  implies  Lemma~\ref{LDP5}(ii)  and  so  it  suffice  to  prove  that  for  a  probability  measure  of  the  form  $\pi=e^{g}\sfrac{\kappa}{\ell}\omega\otimes\omega$,  where  the  Kullback  action  ${\mathfrak  H}_{\kappa/\ell}(\pi\,\|\,\omega)=\langle  g,\,\pi\rangle+ \langle (1-e^{g}),\,\sfrac{\kappa}{\ell}\omega\otimes\omega\rangle$  is  finite.   Fix  any  number $\eps>0$  and  any  neighbourhood  $B_{\pi}\subset\skril(\Omega\times\Omega).$  We  define  the  sequence  of  sets

$$\tilde{\skrig}([n],\,\Omega):=\Big\{y\in \skrig([n],\,\Omega): L_y^{2}\in B_{\pi}\,,\Big|\langle\, g,\,L^2_y\rangle-\langle\, g,\,\pi\rangle\Big|\le  \sfrac{\eps}{2}\Big\}.$$

Observe  that,  for  all $z\in\skrig([n],\,\Omega)$    we  have

$$ \frac{dP_\omega(z)}{d\tilde{P}_\omega(z)}= e^{-n\langle
	\sfrac{1}{2}L^2,\, \tilde{g}\rangle-n\langle\sfrac{1}{2}L^1\otimes
	L^1,\, \tilde{h}_n\rangle+\langle \sfrac{1}{2}L_{\Delta}^{1},\,
	\tilde{h}_n\rangle}>e^{-n\langle\sfrac{1}{2}\pi,\,\log\sfrac{\pi}{\sfrac{\kappa}{\ell}\omega\otimes\omega}\rangle-n\sfrac{1}{2}\langle\sfrac{\kappa}{\ell}\omega\otimes
	\omega,\, (1-\sfrac{\pi}{\sfrac{\kappa}{\ell}\omega\otimes\omega})\rangle}$$

This  gives  $$ \begin{aligned}
P_\omega\Big(\tilde{\skrig}([n],\Omega)\Big)=\int_{\tilde{\skrig}([n],\Omega)}dP_\omega(z)&\ge  \int_{\tilde{\skrig}([n],\,\Omega)}e^{-n\langle\sfrac{1}{2}\pi,\,\log\sfrac{\pi}{\sfrac{\kappa}{\ell}\omega\otimes\omega}\rangle-n\sfrac{1}{2}\langle\sfrac{\kappa}{\ell}\omega\otimes
	\omega,\, (1-\sfrac{\pi}{\sfrac{\kappa}{\ell}\omega\otimes\omega})\rangle+\sfrac{\eps}{2}}d\tilde{P}_\omega(z)\\
&=e^{-n{\mathfrak  H}_{\kappa/\ell}(\pi\,\|\,\omega)/2+n\eps}\tilde{P}_{\omega}\Big(\tilde{\skrig}
([n],\,\Omega)\Big).
\end{aligned}$$
Using  the  law  of  large  numbers  we  have  $\lim_{n\to\infty}\tilde{P}_{\omega}(\tilde{\skrig}([n],\Omega))=1$  which  completes  the  proof.


The  proof  of Theorem~\ref{LDP4}   below,  follows  from  Lemma~\ref{LDP5} above  using  similar  arguments  as  in \cite[p. 544]{BIV15}.
\subsection{Proof  of  Theorem~\ref{LDP4}.}
\begin{proof}
	Note  that  the  empirical  link  measure  is  a  finite  measure  and  so belongs  to  some ball  in  $\skrib_*(\Omega\times\Omega).$ Hence,  without  loss  of  generality  we  may  assume  that  the  set  $\Omega$  in  Theorem~\ref{LDP4}(ii)  is  relatively  compact. See Lemma~~\ref{LDP5} (iii). Choose  any  $\eps>0.$   Then  for  every  functional  $\pi\in \skrim(\Omega\times\Omega)$  we can  find  a  weak  neighbourhood  such  that the  estimate  of  Lemma~\ref{LDP5}(i)  holds.  We  choose  from  all  these  neighbourhood  a  finite  cover  of  $F$ 
	 and  sum up over  the  estimate  in  Theorem~\ref{LDP5}(i)   to  obtain
	
	$$ \lim_{n\to\infty}\frac{1}{n}\log P_n\Big\{\, L_z^2\in F \Big\}\le-\inf_{\pi\in F} I_{1}(\omega,\pi)+\eps.$$
	As $\eps$ was  arbitrarily  chosen  and  the   lower  bound  in  Lemma~\ref{LDP5}(ii) implies  the  lower  bound  in  Theorem~\ref{LDP4}(ii) we  have  the  desired  results  which  ends  the  proof of  the Theorem.
	
\end{proof}

\section{Proof  of Theorem~\ref{LDP2},Proof  of Theorem~\ref{LDP1}  and   Corollary}
\subsection{Proof  of Theorem~\ref{LDP2}}.
Note  that  $dP_n(\omega,\,\pi):=dP_{\omega}(\pi)dP(\omega)$ and the  sequence of   $P_n$ probability  measures  are  exponentially  tight.  See  \cite{DM10} .  Moreover  the  function  $I_1, I$  are lower  semicontinuous  on  the  space  of  measures  $\skrim(\Omega\times\Omega)$.  Therefore,  by  the  Sanov  Theorem,  see  \cite{DZ98}, Theorem~\ref{LDP4}   and  the Theorem for  mixing, see \cite[Theorem~]{Bi04}  we  $(L^1,L^2)$  satisfies an  LDP  with  convex  rate  function  $I$  which  completes  the  proof.
\subsection{Proof  of Theorem~\ref{LDP1}}.
Similarly, we  note  that  the  law  of  the  pair  could  be  written  as  $dP^{n)}(\pi,\omega)=dP_{(\mu,\pi)}(\omega)dP_n(\mu,\pi)$  and  also  the  probability  measures  $P^{(n)}$  are  exponentially  tight.  See \cite{DM10}.  Furthermore,  the  function  $J^1(\pi,\omega)$   and  $J_{(\mu,\pi)}(\omega)$  are  lower  semi-continuous  on  the  space  of  measures  on  $\skrim(\Omega\times\Omega)\times\skrim(\Omega\times\skrin(\Omega))$. Therefore,  by  Theorem~\ref{LDP2}  and  the Theorem for  mixing  large  deviation principles,see \cite[Theorem]{Bi04},  we  have  that  $(L^2, M)$  obeys  an  LDP  with  convex  good  rate  function  $J^1$.  This  ends  the  proof  of  the  Theorem.

\subsection{Corollary~\ref{LDP2}}

{\bf Corollary~\ref{LDP2}(i).}
The  proof  of  Corollary~\ref{LDP2}(i)  follows  from  Theorem~\ref{LDP1}  if  we  take  $p_n\to \kappa/\ell$  and  note  from  the  Doku-Amponsah~\cite[P.p.~]{DM10} that  $\lambda(d)=\{J^1(\pi,\mu):\, \langle d\rangle \le \kappa/\ell,\, \mu=d\}$ which  produces  the  rate  function  in  Corollary~\ref{LDP2}(i).\\

{\bf Corollary~\ref{LDP2}(ii).} Similarly, the  proof  of  the (ii)  part  follows  from  the   Corollary~\ref{LDP2}(i) by  applying  the  contraction  principle  to  the  linear  mapping  $U(d)=d(0)$ and  solve  the  variational  problem  $\xi(y)=\inf\{\lambda(d): d(0)=y\}$  as  in  Doku-Amponsah~\cite[pp. ]{DA16} to  obtain  the  desired  form  of  the  rate  function  in  Corollary~\ref{LDP2}(ii).

This  completes  the  proof  of  the  Corollary.


{\bf \Large Conflict  of  Interest}

The  author  declares  that  he has  no  conflict  of  interest.\\





\end{document}